\documentclass[a4paper,12pt,leqno]{article}
\usepackage{amsmath,amssymb,latexsym,amscd,amsthm,amsxtra,graphicx, fancybox}

\begin{document}


\newtheorem{thm}{Theorem}[subsection]
\newtheorem{lem}[thm]{Lemma}
\newtheorem{cor}[thm]{Corollary}
\newtheorem{prop}[thm]{Proposition}
\newtheorem{prob}[thm]{Problem}
\newtheorem{clm}{Claim}
\newtheorem{dfn}[thm]{Definition}
\newtheorem{conj}[thm]{Conjecture}
\newtheorem{wh}[thm]{Working hypothesis}


\newcommand{\pf}{{\it Proof} \,\,\,\,\,\,}
\newcommand{\remark}{{\it Remark}. \,\,}

\newcommand{\eqb}{\begin{equation}}
\newcommand{\eqe}{\end{equation}}

\newcommand{\cpx}{{\mathbb C}}
\newcommand{\rel}{{\mathbb R}}
\newcommand{\rat}{{\mathbb Q}}
\newcommand{\itg}{{\mathbb Z}}
\newcommand{\nat}{{\mathbb N}}
\newcommand{\qrt}{{\mathbb H}}
\newcommand{\qtr}{{\mathbb H}}
\newcommand{\dbP}{{\mathbb P}}
\newcommand{\dbS}{{\mathbb S}}
\newcommand{\dbO}{{\mathbb O}}
\newcommand{\dbK}{{\mathbb K}}

\newcommand{\mathbbP}{{\mathbb P}}
\newcommand{\mathbbA}{{\mathbb A}}
\newcommand{\mathbbB}{{\mathbb B}}
\newcommand{\mathbbC}{{\mathbb C}}
\newcommand{\mathbbD}{{\mathbb D}}
\newcommand{\mathbbE}{{\mathbb E}}
\newcommand{\mathbbF}{{\mathbb F}}
\newcommand{\mathbbG}{{\mathbb G}}
\newcommand{\mathbbH}{{\mathbb H}}
\newcommand{\mathbbI}{{\mathbb I}}
\newcommand{\mathbbJ}{{\mathbb J}}
\newcommand{\mathbbK}{{\mathbb K}}
\newcommand{\mathbbL}{{\mathbb L}}
\newcommand{\mathbbM}{{\mathbb M}}
\newcommand{\mathbbN}{{\mathbb N}}
\newcommand{\mathbbO}{{\mathbb O}}
\newcommand{\mathbbQ}{{\mathbb Q}}
\newcommand{\mathbbR}{{\mathbb R}}
\newcommand{\mathbbS}{{\mathbb S}}
\newcommand{\mathbbT}{{\mathbb T}}
\newcommand{\mathbbU}{{\mathbb U}}
\newcommand{\mathbbV}{{\mathbb V}}
\newcommand{\mathbbW}{{\mathbb W}}
\newcommand{\mathbbX}{{\mathbb X}}
\newcommand{\mathbbY}{{\mathbb Y}}
\newcommand{\mathbbZ}{{\mathbb Z}}

\newcommand{\arn}{\Vec{n}}

\newcommand{\uDl}{\underline{\Delta}}
\newcommand{\uPi}{\underline{\Pi}}
\newcommand{\oDl}{\overline{\Delta}}
\newcommand{\uW}{\underline{W}}
\newcommand{\oW}{\overline{W}}
\newcommand{\uqca}{\underline{\cal Q}}
\newcommand{\upca}{\underline{\cal P}}
\newcommand{\oqca}{\overline{\cal Q}}
\newcommand{\opca}{\overline{\cal P}}
\newcommand{\urho}{\underline{\rho}}
\newcommand{\orho}{\overline{\rho}}
\newcommand{\ubbb}{\underline{\mathfrak b}}
\newcommand{\uqqq}{\underline{\mathfrak q}}
\newcommand{\uppp}{\underline{\mathfrak p}}
\newcommand{\uvvv}{\underline{\mathfrak v}}
\newcommand{\unnn}{\underline{\mathfrak n}}
\newcommand{\obbb}{\overline{\mathfrak b}}
\newcommand{\onnn}{\overline{\mathfrak n}}
\newcommand{\ola}{\overline{\lambda}}
\newcommand{\dola}{\overline{\overline{\lambda}}}
\newcommand{\ula}{\underline{\lambda}}
\newcommand{\dula}{\underline{\underline{\lambda}}}
\newcommand{\hla}{{\widehat{\lambda}}}

\newcommand{\bGL}{{\mathbb GL}}

\newcommand{\bfa}{\mbox{\underline{\bf a}}}
\newcommand{\bfb}{\mbox{\underline{\bf b}}}
\newcommand{\bfc}{\mbox{\underline{\bf c}}}
\newcommand{\bfd}{\mbox{\underline{\bf d}}}
\newcommand{\bfe}{\mbox{\underline{\bf e}}}
\newcommand{\bff}{\mbox{\underline{\bf f}}}
\newcommand{\bfg}{\mbox{\underline{\bf g}}}
\newcommand{\bfh}{\mbox{\underline{\bf h}}}
\newcommand{\bfi}{\mbox{\underline{\bf i}}}
\newcommand{\bfj}{\mbox{\underline{\bf j}}}
\newcommand{\bfk}{\mbox{\underline{\bf k}}}
\newcommand{\bfl}{\mbox{\underline{\bf l}}}
\newcommand{\bfm}{\mbox{\underline{\bf m}}}
\newcommand{\bfn}{\mbox{\underline{\bf n}}}
\newcommand{\bfo}{\mbox{\underline{\bf o}}}
\newcommand{\bfp}{\mbox{\underline{\bf p}}}
\newcommand{\bfq}{\mbox{\underline{\bf q}}}
\newcommand{\bfr}{\mbox{\underline{\bf r}}}
\newcommand{\bfs}{\mbox{\underline{\bf s}}}
\newcommand{\bft}{\mbox{\underline{\bf t}}}
\newcommand{\bfu}{\mbox{\underline{\bf u}}}
\newcommand{\bfv}{\mbox{\underline{\bf v}}}
\newcommand{\bfw}{\mbox{\underline{\bf w}}}
\newcommand{\bfx}{\mbox{\underline{\bf x}}}
\newcommand{\bfy}{\mbox{\underline{\bf y}}}
\newcommand{\bfz}{\mbox{\underline{\bf z}}}
\newcommand{\bfmu}{{\underline{\bf \mu}}}

\newcommand{\bpa}{\mbox{\underline{\bf a}}}
\newcommand{\bpb}{\mbox{\underline{\bf b}}}
\newcommand{\bpc}{\mbox{\underline{\bf c}}}
\newcommand{\bpd}{\mbox{\underline{\bf d}}}
\newcommand{\bpe}{\mbox{\underline{\bf e}}}
\newcommand{\bpf}{\mbox{\underline{\bf f}}}
\newcommand{\bpg}{\mbox{\underline{\bf g}}}
\newcommand{\bph}{\mbox{\underline{\bf h}}}
\newcommand{\bpi}{\mbox{\underline{\bf i}}}
\newcommand{\bpj}{\mbox{\underline{\bf j}}}
\newcommand{\bpk}{\mbox{\underline{\bf k}}}
\newcommand{\bpl}{\mbox{\underline{\bf l}}}
\newcommand{\bpm}{\mbox{\underline{\bf m}}}
\newcommand{\bpn}{\mbox{\underline{\bf n}}}
\newcommand{\bpo}{\mbox{\underline{\bf o}}}
\newcommand{\bbpp}{\mbox{\underline{\bf p}}}
\newcommand{\bpq}{\mbox{\underline{\bf q}}}
\newcommand{\bpr}{{\mbox{\underline{\bf r}}}}
\newcommand{\bps}{\mbox{\underline{\bf s}}}
\newcommand{\bpt}{\mbox{\underline{\bf t}}}
\newcommand{\bpu}{\mbox{\underline{\bf u}}}
\newcommand{\bpv}{\mbox{\underline{\bf v}}}
\newcommand{\bpw}{\mbox{\underline{\bf w}}}
\newcommand{\bpx}{\mbox{\underline{\bf x}}}
\newcommand{\bpy}{\mbox{\underline{\bf y}}}
\newcommand{\bpz}{\mbox{\underline{\bf z}}}

\newcommand{\eps}{\varepsilon}

\newcommand{\T}{\Theta}
\newcommand{\Th}{\Theta}
\newcommand{\Ta}{\Theta^\alpha}

\newcommand{\zzz}{{\mathfrak z}}
\newcommand{\fff}{{\mathfrak f}}
\newcommand{\www}{{\mathfrak w}}
\newcommand{\vvv}{{\mathfrak v}}
\renewcommand{\ggg}{{\mathfrak g}}
\newcommand{\gggg}{{\mathfrak g}}
\newcommand{\kkk}{{\mathfrak k}}
\newcommand{\aaa}{{\mathfrak a}}
\newcommand{\saa}{\mbox{}^s{\mathfrak a}}
\newcommand{\aas}{{\mathfrak a}_{\T}}
\newcommand{\aasa}{{\mathfrak a}_{\T^\alpha}}
\newcommand{\bbb}{{\mathfrak b}}
\newcommand{\sbb}{\mbox{}^s{\mathfrak p}}
\newcommand{\spp}{\mbox{}^s{\mathfrak p}}
\newcommand{\llll}{{\mathfrak l}}
\newcommand{\llls}{{\mathfrak l}_{\T}}
\newcommand{\rrr}{{\mathfrak r}}
\newcommand{\tll}{\widetilde{\mathfrak l}}
\newcommand{\eee}{{\mathfrak e}}
\newcommand{\ccc}{{\mathfrak c}}
\newcommand{\ttt}{{\mathfrak t}}
\newcommand{\stt}{\mbox{}^s{\mathfrak t}}
\newcommand{\lls}{{\mathfrak l}_{\T}}
\newcommand{\obb}{\bar{\mathfrak b}}
\newcommand{\aad}{{\mathfrak a}^\ast_{\T}}
\newcommand{\ads}{{\mathfrak a}^\ast_{\T}}
\newcommand{\nnn}{{\mathfrak n}}
\newcommand{\nns}{{\mathfrak n}_{\T}}
\newcommand{\uuu}{{\mathfrak u}}
\newcommand{\tuu}{\widetilde{\mathfrak u}}
\newcommand{\buu}{\bar{\mathfrak u}}
\newcommand{\bnn}{\bar{\mathfrak n}}
\newcommand{\bab}{\bar{\mathfrak b}}
\newcommand{\bns}{\bar{\mathfrak n}_{\T}}
\newcommand{\mmm}{{\mathfrak m}}
\newcommand{\smm}{\mbox{}^s{\mathfrak m}}
\newcommand{\mms}{{\mathfrak m}_{\T}}
\newcommand{\jjj}{{\mathfrak j}}
\newcommand{\hhh}{{\mathfrak h}}
\newcommand{\shh}{\mbox{}^s{\mathfrak h}}
\newcommand{\uhh}{\mbox{}^u{\mathfrak h}}
\newcommand{\hhd}{{\mathfrak h}^\ast}
\newcommand{\qqq}{{\mathfrak q}}
\newcommand{\ppp}{{\mathfrak p}}
\newcommand{\tpp}{\widetilde{\mathfrak p}}
\newcommand{\tnn}{\widetilde{\mathfrak n}}
\newcommand{\tqq}{\widetilde{\mathfrak q}}
\newcommand{\pps}{{\mathfrak p}_{\T}}
\newcommand{\bapp}{\bar{\mathfrak p}}
\newcommand{\baps}{\bar{\mathfrak p}_{\T}}
\newcommand{\sss}{{\mathfrak s}}
\newcommand{\ooo}{{\mathfrak o}}
\newcommand{\ddd}{{\mathfrak d}}
\newcommand{\Sth}{{\mathfrak S}}

\newcommand{\ggc}{{\mathfrak g}_{\sf c}}
\newcommand{\bbc}{{\mathfrak b}_{\sf c}}
\newcommand{\llc}{{\mathfrak l}_{\sf c}}
\newcommand{\uuc}{{\mathfrak u}_{\sf c}}
\newcommand{\nsc}{{({\mathfrak n}_{\T})}_{\sf c}}
\newcommand{\nnc}{{({\mathfrak n}_{\T})}_{\sf c}}
\newcommand{\buc}{\bar{\mathfrak u}_{\sf c}}
\newcommand{\bnc}{{(\bar{\mathfrak n}_{\T})}_{\sf c}}
\newcommand{\hhc}{{\mathfrak h}_{\sf c}}
\newcommand{\hcd}{{{\mathfrak h}_{\sf c}}^\ast}
\newcommand{\kkc}{{\mathfrak k}_{\sf c}}
\newcommand{\ppc}{{\mathfrak p}_{\sf c}}
\newcommand{\pcs}{({\mathfrak p}_{\T})_{\sf c}}
\newcommand{\ncs}{({\mathfrak n}_{\T})_{\sf c}}
\newcommand{\bpp}{\bar{\mathfrak p}}
\newcommand{\bqq}{\bar{\mathfrak q}}

\newcommand{\spl}{\mbox{}^s}

\newcommand{\gl}{{\mathfrak g}{\mathfrak l}}
\newcommand{\gll}{{\mathfrak gl}_L}
\newcommand{\glr}{{\mathfrak gl}_R}
\newcommand{\so}{\sss\ooo}

\newcommand{\sA}{{\mbox{}^s{A}}}
\newcommand{\sM}{{\mbox{}^s{M}}}
\newcommand{\sT}{{\mbox{}^s{T}}}
\newcommand{\sH}{{\mbox{}^s{H}}}
\newcommand{\sB}{{\mbox{}^s{P}}}
\newcommand{\sP}{{\mbox{}^s{P}}}
\newcommand{\uH}{{\mbox{}^u{H}}}

\newcommand{\bi}{{\bold i}}
\newcommand{\fsp}{{\mathfrak s}{\mathfrak p}}

\newcommand{\Gc}{{G_{\mathbb C}}}
\newcommand{\Gac}{{G_{\mathbb C}^{\rm ad}}}
\newcommand{\Gca}{{G_{\mathbb C}^{\rm ad}}}
\newcommand{\Gf}{{G^{\flat}}}
\newcommand{\Gs}{{G^{\sharp}}}
\newcommand{\Ga}{{G^{\ast}}}
\newcommand{\Kc}{{K_{\mathbb C}}}
\newcommand{\Kcf}{{K_{\mathbb C}^\flat}}
\newcommand{\Kf}{{K^\flat}}
\newcommand{\Ks}{{K^\sharp}}
\newcommand{\Ka}{{K^\ast}}
\newcommand{\Kcs}{{K_{\mathbb C}^\sharp}}
\newcommand{\Kca}{{K_{\mathbb C}^\ast}}
\newcommand{\Bc}{{B_{\mathbb C}}}
\newcommand{\bBc}{{\bar{B}_{\sf C}}}
\newcommand{\Pc}{{P_{\mathbb C}}}
\newcommand{\bPc}{{\bar{P}_{\mathbb C}}}
\newcommand{\Nc}{{N_{\mathbb C}}}
\newcommand{\bNc}{{\bar{N}_{\mathbb C}}}

\newcommand{\hol}{{\cal O}}
\newcommand{\dif}{{\cal D}}
\newcommand{\ana}{{\cal A}}

\newcommand{\tpca}{\widetilde{\cal P}}
\newcommand{\aca}{{\cal A}}
\newcommand{\bca}{{\cal B}}
\newcommand{\cca}{{\cal C}}
\newcommand{\dca}{{\cal D}}
\newcommand{\diff}{{\cal D}}
\newcommand{\eca}{{\cal E}}
\newcommand{\fca}{{\cal F}}
\newcommand{\gca}{{\cal G}}
\newcommand{\hca}{{\cal H}}
\newcommand{\ica}{{\cal I}}
\newcommand{\jca}{{\cal J}}
\newcommand{\kca}{{\cal K}}
\newcommand{\lca}{{\cal L}}
\newcommand{\mca}{{\cal M}}
\newcommand{\nca}{{\cal N}}
\newcommand{\oca}{{\cal O}}
\newcommand{\pca}{{\cal P}}
\newcommand{\qca}{{\cal Q}}
\newcommand{\rca}{{\cal R}}
\newcommand{\sca}{{\cal S}}
\newcommand{\tca}{{\cal T}}
\newcommand{\uca}{{\cal U}}
\newcommand{\vca}{{\cal V}}
\newcommand{\wca}{{\cal W}}
\newcommand{\xca}{{\cal X}}
\newcommand{\yca}{{\cal Y}}
\newcommand{\zca}{{\cal Z}}

\newcommand{\HG}{{\cal H}_G}
\newcommand{\PR}{{\cal P}r}
\newcommand{\gS}{{\mathfrak S}}
\newcommand{\ii}{\sqrt{-1}}
\newcommand{\real}{\mbox{Re}}
\newcommand{\Real}{\mbox{Re}}
\newcommand{\res}{\mbox{res}}
\newcommand{\supp}{\mbox{supp}}
\newcommand{\rank}{\mbox{rank}}
\newcommand{\card}{\mbox{card}}
\newcommand{\Ad}{\mbox{Ad}}
\newcommand{\ad}{\mbox{ad}}
\newcommand{\dg}{{\deg}}
\newcommand{\Hom}{\mbox{Hom}}
\newcommand{\End}{\mbox{End}}
\newcommand{\Ext}{\mbox{Ext}}
\newcommand{\fgt}{\mbox{Fgt}}
\newcommand{\JH}{\mbox{JH}}
\newcommand{\Dim}{\mbox{Dim}}
\newcommand{\gr}{\mbox{gr}}
\newcommand{\Ass}{\mbox{Ass}}
\newcommand{\WF}{\mbox{WF}}
\newcommand{\AS}{\mbox{AS}}
\newcommand{\Ann}{\mbox{Ann}}
\newcommand{\LAnn}{\mbox{LAnn}}
\newcommand{\RAnn}{\mbox{RAnn}}
\newcommand{\tr}{\mbox{tr}}
\newcommand{\hht}{\mbox{ht}}
\newcommand{\Mod}{{\sf Mod}}
\newcommand{\sgn}{\mbox{sgn}}
\newcommand{\pro}{\mbox{pro}}
\newcommand{\Ind}{\mbox{\sf Ind}}
\newcommand{\SO}{\mbox{SO}}
\newcommand{\Oo}{\mbox{O}}
\newcommand{\SOo}{{\mbox{SO}_0}}
\newcommand{\GL}{\mbox{GL}}
\newcommand{\Spp}{\mbox{Sp}}
\newcommand{\U}{\mbox{U}}
\newcommand{\triv}{{\mbox{triv}}}
\newcommand{\Rea}{{\mbox{Re}}}

\newcommand{\Res}{\mbox{Res}}
\newcommand{\diag}{\mbox{diag}}

\newcommand{\PP}{{\sf P}}
\newcommand{\PS}{{\sf P}_{\T}^{++}}
\newcommand{\PSP}{{\sf P}_{\T^\prime}^{++}}
\newcommand{\PPS}{{\sf P}_{\T}^{++}}
\newcommand{\PPSA}{{\sf P}_{\T^\alpha}^{++}}
\newcommand{\PPSK}{{\sf P}_{\T^k}^{++}}

\newcommand{\Q}{{\sf Q}}

\newcommand{\leqs}{\leqslant}
\newcommand{\geqs}{\geqslant}
\title{{\sf Homomorphisms between scalar generalized Verma modules of ${\mathfrak gl} (n,\cpx)$}}
\author{{\sf  Hisayosi Matumoto}
\\Graduate School of Mathematical Sciences\\ University of Tokyo\\ 3-8-1
Komaba, Tokyo\\ 153-8914, JAPAN\\ e-mail: hisayosi@ms.u-tokyo.ac.jp}
\date{}
\maketitle
\begin{abstract}
In this article, we classify  the homomorphisms between scalar generalized
Verma modules of $\gl(n,\cpx)$.  In fact such homomorphisms are compositions of elementary homomorphisms.
\footnote{This work was supported by Grant-in-Aid for Scientific Research (No.\ 2054001, No.\ 260006)}
\end{abstract}

\setcounter{section}{0}
\section*{\S\,\, 0.\,\,\,\, Introduction}
\setcounter{subsection}{0}

An induced module of a complex reductive Lie algebra from a  one-dimensional representation of a parabolic subalgebra is called a scalar generalized Verma module.
In this article, we give a classification of homomorphisms between scalar generalized Verma modules of $\gl(n,\cpx)$.

In \cite{[Vr]}, Verma constructed  homomorphisms between
Verma modules of a complex reductive Lie algebra associated with root reflections.
Bernstein, I.\ M.\ Gelfand, and S.\ I.\ Gelfand proved that all the nontrivial homomorphisms between Verma modules are compositions of homomorphisms constructed by Verma.  (\cite{[BGG]}) 

Later, Lepowsky studied generalized Verma modules.
In particular, Lepowsky (\cite{[L]}) constructed a class of homomorphisms between scalar generalized Verma modules 
associated to the parabolic
subalgebras which are the complexifications of the minimal parabolic
subalgebras of real reductive Lie algebras.  They correspond to reflections with respect to the restricted roots.

In \cite{[M2]}, elementary homomorphisms (see Proposition 2.3.1 below)
 between scalar generalized Verma modules
are introduced.  They can be regarded as a generalization of homomorphisms introduced by Verma and Lepowsky.
The main theorem of this article is as follows. 

\medskip

{\bf Theorem} (Theorem 2.3.2) \,\,\,\,
{\it Non-zero homomorphisms between scalar generalized Verma modules of $\gl(n, \cpx)$ are compositions of elementary homomorphisms.}
\medskip

This result confirms Conjecture A in \cite{[M3]} for $\gl(n,\cpx)$.

The main ingredient of our proof of the theorem is the translation principle in mediocre regions studied by Vogan, Kobayashi, and Trapa (\cite{[K]}, \cite{[T]}, \cite{[VU]}, \cite{[VI]}, \cite{[Vd]}).
A key result is non-existence of  certain homomorphisms.  We assume existence of such homomorphisms.  Then, applying translation functors, we obtain  homomorphisms at very degenerated parameter, of which we easily see the non-existence.
We also use a result of Borho-Jantzen (\cite{[BJ]} 5.10) in order to show that a non-zero homomorphism exists between scalar generalized Verma modules only if their parameters are in the same $W(\T)$-orbit. (See Lemma 2.4.1 below.)

\setcounter{section}{1}
\setcounter{subsection}{0}

\section*{\S\,\, 1.\,\,\,\,Notations and Preliminaries }

\subsection{General notations}

In this article, we use the following notations and conventions.

As usual we denote the complex number field, the real number field, the
ring of (rational) integers, and the set of non-negative integers by
$\cpx$, $\rel$, $\itg$, and $\nat$ respectively.
We denote by $\emptyset$ the empty set.
For any (non-commutative) $\cpx$-algebra $R$,  ``$R$-module'' means ``left $R$-module'', and sometimes we denote by $0$ (resp.\ $1$) the trivial $R$-module $\{0\}$ (resp.\ $\cpx$).
We denote by $\Ann_R(M)$ the annihilator of $M$ in $R$.
Often, we identify a (small) category and the set of its objects.
Hereafter ``$\dim$'' means the dimension as a complex vector space, and ``$\otimes$'' (resp. $\Hom$) means the tensor product over $\cpx$ (resp. the space of $\cpx$-linear mappings), unless we specify otherwise.
For a complex vector space $V$, we denote by $V^\ast$ the dual vector space.
For $a,b\in\cpx$, ``$a\leqslant b$'' means that $a,b\in\rel$ and $a\leqslant b$.
We denote by $A- B$ the set theoretical difference.
$\card A$ means the cardinality of a set $A$.
We denote by $\delta_{i,j}$ the Kronecker delta.  Namely,
\begin{equation}\delta_{i.j}=\left\{ \begin{array}{l} 1 \,\,\,\, \hbox{(if $i=j$)} \\ 0 \,\,\,\, \hbox{(if $i\neq j$)}\end{array}\right.\end{equation}

Let $\qqq$ be a Lie algebra and let $\rrr$ be a subalgebra of ${\mathfrak q}$. For a ${\mathfrak q}$-module $V$, we denote by $V|_{\mathfrak r}$ the restriction to $\rrr$.  For a Lie algebra $\qqq$, we denote by $U(\qqq)$ its universal enveloping algebra.

\subsection{Notations on $\ggg\llll(n,\cpx)$}
In this article, we fix a positive integer $n$, which is greater than $1$.  
We put $\ggg=\ggg\llll(n,\cpx)$.
Let $\hhh$ be the Cartan subalgebra of $\ggg$ consisting of the diagonal matrices and let $\bbb$ be the Borel subalgebra of $\ggg$ consisting of the upper triangular matrices. We denote by $\Delta$ the root system with respect to $(\ggg,\hhh)$. We choose positive root system $\Delta^+$ corresponding to $\bbb$. 
Let $W$ be the Weyl group of the pair $(\ggg, \hhh)$.  We denote by $s_\alpha$ the refection with respect to $\alpha\in\Delta$.   
Let $\langle\,\,,\,\,\rangle$ be a non-degenerate invariant bilinear form on $\ggg$ defined as follows.
$$\langle X,Y\rangle=\tr XY     \,\,\,\,\,\,  (X,Y\in \ggg).$$
This bilinear form induces a $W$-invarinant non-degenerate bilinear form on $\hhh^\ast$ in a usual way.  We also denote it by the same symbol $\langle\,\,,\,\,\rangle$.

For $1\leqslant i,j\leqslant n$, we denote by $E_{i,j}$ the matrix element for $(i,j)$.  Namely, $E_{i,j}$ is the $n\times n$ matrix whose $p,q$-entry is $\delta_{p,i}\delta_{q.j}$.  So, $\{E_{i,j}\mid 1\leqslant i,j\leqslant n\}$ is a basis of $\ggg$ as a $\cpx$-vector space.
 For $\alpha\in \Delta$, we denote by $\ggg_\alpha$ the corresponding root space.  So, we have $\ggg_{e_i-e_j}=\cpx E_{i,j}$ for all $1\leqslant i,j\leqslant n$ such that $i\neq j$.
For $1\leqslant i\leqslant n$, we put $H_i=E_{i,i}$.  Then, $\{H_1,...,H_n\}$ is a basis of $\hhh$.
We denote by  $e_1,...,e_n$ the dual basis of $\hhh^\ast$ corresponding to a basis $H_1,....,H_n$.
Then, $e_1,...,e_n$ form an orthonormal basis with respect to  $\langle\,\,,\,\,\rangle$ and we see that 
\begin{align*}
\Delta^+ = \{ e_i-e_j\mid 1\leqslant i< j\leqslant
n, i\neq j\}.
\end{align*}
We put $\alpha_i=e_i-e_{i+1}$  \,\,\, $(1\leqslant i <n)$.
We see that the basis of $\Delta$ with respect to $\Delta^+$ is $\Pi=\{\alpha_1,...,\alpha_{n-1}\}$.
We identify the Weyl group $W$ with the $n$-th symmetric group ${\mathfrak S}_n$ via $\sigma e_i=e_{\sigma(i)}$  \,\,\,\,($1\leqslant i\leqslant n$).

\subsection{Notations on parabolic subalgebras}

We fix $\T\subsetneq\Pi$ and write $\Pi-\T=\{\alpha_{i_1},...,\alpha_{i_{k-1}}\}$, where $1\leqslant i_1<\cdots<i_{k-1}\leqslant n-1$.
We put $i_0=0$ and $i_k=n$.
For $1\leqslant j\leqslant k$, we put $n_j=i_{j-1}-i_j$, so that we have $n=n_1+\cdots+n_k$.

Let $\langle\T\rangle$ be the set of the elements of $\Delta$ which are written by linear combinations of elements of $\T$ over $\itg$.
Put
$\aas = \{H\in\hhh\mid \forall \alpha\in \T \,\,\,\alpha(H)=0\}$,
$\lls = \hhh +\sum_{\alpha\in\langle{\T}\rangle}\ggg_\alpha$,  
$\nns = \sum_{\alpha\in\Delta^+-\langle{\T}\rangle}\ggg_\alpha$,  
$\pps = \lls +\nns$.
Then $\pps$ is a parabolic subalgebra of $\ggg$ which contains $\bbb$.
Conversely, for an arbitrary parabolic subalgebra $\ppp\supseteq\bbb$, there exists some $\T \subseteq \Pi$ such that $\ppp =\pps$.
$\lls$ (resp.\ $\pps$)  consists of the diagonal (resp.\ upper-triangular) block $n\times n$ matrices with respect to a composition $n=n_1+\cdots+n_k$.  So, we see 
$$ \lls\cong \ggg\llll(n_1,\cpx)\oplus\cdots\oplus\ggg\llll(n_k,\cpx).$$

We denote by $W_\T$ the Weyl group for $(\lls,\hhh)$. 
$W_\T$ is identified with a subgroup of $W$ generated by $\{s_\alpha\mid \alpha\in\T\}$.
We denote by $w_\T$ the longest element of $W_\T$.
Using the invariant non-degenerate bilinear form $\langle\,\,,\,\,\rangle$, we regard  ${\aaa_{\T}}^\ast$ as a subspace of $\hhh^\ast$. 

Put
$\rho_{\T} = \frac{1}{2}(\rho-w_{\T}\rho)$. 


Define
\begin{align*}
\PS & = \{\lambda\in\hhd\mid \forall\alpha\in \T\,\,\,\,\,
 \langle\lambda,{\alpha}^\vee\rangle\in\{1,2,...\}\}\\
{}^\circ\PS & =\{\lambda\in\hhd\mid \forall\alpha\in \T\,\,\,\,\,
 \langle\lambda,{\alpha}^\vee\rangle=1\}
\end{align*}
We easily have
\begin{align*}
{}^\circ\PS =\{\rho_{\T}+\mu\mid \mu\in\ads\}.
\end{align*}
For $\mu\in\hhd$ such that $\mu+\rho\in\PS$, we denote by $\sigma_\T(\mu)$ the irreducible finite-dimensional $\lls$-representation whose highest weight is $\mu$.
Let $E_\T(\mu)$ be the representation space of $\sigma_\T(\mu)$.
We define a left action of $\nns$ on $E_\T(\mu)$ by $X\cdot v =0$ for all $X\in\nns$ and $v\in E_\T(\mu)$.
So, we regard $E_\T(\mu)$ as a $U(\pps)$-module.

For $\mu\in\PS$, we define a generalized Verma module (\cite{[L4]}) as follows.
\[ M_\T(\mu) = U(\gggg)\otimes_{U(\pps)}E_\T({\mu-\rho}).\]

We see that $\dim E_{\T}(\mu-\rho)=1$ if and only if $\mu\in {}^\circ\PS$.
If $\mu\in {}^\circ\PS$, we call $M_\T(\mu)$ a scalar generalized
Verma module (\cite{[Boe]}).

For $1\leqslant i\leqslant k$, we put $n^\ast_i=n_1+\cdots+n_i$ and put $n_0^\ast=0$.
For $\lambda_1,...,\lambda_k\in \cpx$, we put 
$$ [\lambda_1,...,\lambda_k]=\sum_{j=1}^k\sum_{s=1}^{n_j}(\lambda_j-s+1)e_{n_{j-1}^\ast+s}.$$
Then, we have ${}^\circ\PS=\{  [\lambda_1,...,\lambda_k]\mid  \lambda_1,...,\lambda_k\in\cpx\}.$

\subsection{Notations on infinitesimal characters}
Finally, we fix notations for infinitesimal characters.
We denote by $Z(\ggg)$ the center of $U(\ggg)$.  We denote by $\chi_\lambda$ the image of $\lambda\in\hhd$ under the Harish-Chandra isomorphism from $W\backslash \hhd$ to $\Hom(Z(\ggg), \cpx)$.
It is well-known that $Z(\ggg)$ acts on $M(\lambda)$ by $\chi_\lambda : Z(\ggg) \rightarrow \cpx$ for all $\lambda\in\hhd$.
We denote by $\hbox{\bf Z}_\lambda$ the kernel of $\chi_\lambda$ in $Z(\ggg)$.
Let $M$ be a $U(\ggg)$-module and $\lambda\in\hhh^\ast$.
We say that $M$ has an infinitesimal character $\lambda$ if $Z(\ggg)$ acts on $M$ by $\chi_\lambda$.
We say that $M$ has an generalized infinitesimal character $\lambda$ if there exists  some positive integer $N$ such that $(p-\chi_\lambda(p))^Nv=0$ for all $p\in Z(\ggg)$ and $v\in M$.
We denote by $\mca_\lambda$ the full  subcategory of the category of the $U(\ggg)$-modules consisting of the $U(\ggg)$-modules with generalized infinitesimal character $\lambda$.

For example, a generalized Verma module $M_\T(\mu)$ has an infinitesimal character $\mu$.

A $U(\ggg)$-modules $M$ is called $Z(\ggg)$-finite, if the annihilator $\Ann_{Z(\ggg)}(M)$ of $M$ in $Z(\ggg)$ is finite-codimensional in $Z(\ggg)$.
We denote by $\mca_{Zf}$ the full  subcategory of the category of the $U(\ggg)$-modules consisting of $Z(\ggg)$-finite $U(\ggg)$-modules.
We have a direct sum decomposition of the category.
$$\mca_{Zf}=\bigoplus_{\lambda\in W\backslash \hhd}\mca_\lambda.$$
We denote by ${\bf P}_\lambda : \mca_{Zf}\rightarrow\mca_{\lambda}$ the projection functor with respect to the above direct sum decomposition. ${\bf P}_\lambda$ is obviously an exact functor.

\setcounter{section}{2}
\setcounter{subsection}{0}

\section*{\S\,\, 2.\,\,\,\,Formulation of the main result }

We retain the notation of \S 1.
In particular, $\T$ is a proper subset of $\Pi$.

\subsection{Formulation of the problem}

In \cite{[L3]}, Lepowsky proved that any non-zero  homomorphisms between scalar generalized Verma modules is injective and unique up to a scalar multiplication.  For $\mu,\nu\in{}^\circ \PS$, we write $M_{\T}(\mu)\subseteq M_{\T}(\nu)$ if there is a non-zero homomorphism of $M_{\T}(\mu)$ to $M_{\T}(\nu)$.

The classification problem of homomorphisms between generalized Verma modules is reduced to the following problem.

{\bf Problem 1} \,\, Let $\mu,\nu\in {}^\circ\PS$.
 When is $M_{\T}(\mu)\subseteq M_{\T}(\nu)$ ?

The main result gives a solution to the above problem.
In order to describe it, we explain some notions.

\subsection{Weyl group for $\pps$}

The material in this subsection is more or less  a special case of the results in \cite{[Lu]} and \cite{[H]}.

We put 
\[ W({\T})=\{w\in W\mid w\T= \T \}.\]
For $1\leqslant p<q\leqslant k$ such that $n_p=n_q$, we define $\sigma_{p.q}\in W(={\mathfrak S}_n)$ as follows.
$$ \sigma_{p,q}(j)=
\begin{cases}
i_{q-1}+(j-i_{p-1})  & \hbox{if $i_{p-1}<j\leqslant i_p$} \\
i_{p-1}+(j-i_{q-1})  & \hbox{if $i_{q-1}<j\leqslant i_q$} \\
j & \hbox{otherwise}.
\end{cases}
$$
For $\lambda_1,...\lambda_k\in\cpx$, we have 
$$\sigma_{p,q}[\lambda_1,...,\lambda_k]=[\lambda_1,...,\overset{p}{\widehat{\lambda_q}},...,\overset{q}{\widehat{\lambda_p}},...,\lambda_k].$$

We easily see $\sigma_{p,q}\in W(\T)$.  Moreover, we see that $W(\T)$ is generated by $\sigma_{p,q}$  \,\,\, ($1\leqslant p<q\leqslant k$ and $n_p=n_q$).   We put $I_r=\{j\in\itg\mid n_j=r, 1\leqslant j\leqslant k\}$ and $\Upsilon=\{r\in\itg\mid r\geqslant 1,  I_r\neq\emptyset\}$.
For $r\in\Upsilon$, we denote by $W(\T)_r$ the subgroup of $W(\T)$ generated by $\{\sigma_{p,q}\mid p,q\in I_r, p<q\}$.
We easily see $W(\T)_r\cong {\mathfrak S}_{\card I_r}$ and
\begin{align*}
W(\T)=\prod_{r\in\Upsilon}W(\T)_r. \tag{$\sharp$}
\end{align*}
In other words, $W(\T)$ can be identified as follows.
\begin{align*}
W(\T)\cong\{\sigma\in{\mathfrak S}_k\mid n_{\sigma(i)}=n_i  \,\,\, (1\leqslant i\leqslant k)\}.
\tag{$\flat$}
\end{align*}
For $w\in W(\T)$, we denote by $\bar{w}$ the corresponding element in ${\mathfrak S}_k$.
For example, $\overline{\sigma_{p,q}}$ is the transposition exchanging $p$ and $q$.
We also see:
$$w[\lambda_1,...,\lambda_k]=[\lambda_{\bar{w}^{-1}(1)},...,\lambda_{\bar{w}^{-1}(k)}].$$

\subsection{The main result}

Since the homomorphisms between scalar generalized Verma modules associated with maximal parabolic subalgebras of $\ggg\llll(n,\cpx)$ are classified by Boe (\cite{[Boe]}), the following result follows from {\cite{[M2]} Theorem 5.1.2.

\begin{prop}
Let $\lambda_1,...,\lambda_k\in\cpx$ and $1\leqslant p<q\leqslant k$ be such that $n_p=n_q$ and $\lambda_p-\lambda_q\in\nat$.
Then, we have
$$ M_\T(\sigma_{p,q}[\lambda_1,...\lambda_k])\subseteq M_\T([\lambda_1,...\lambda_k]).$$
\end{prop}

The embedding of a scalar generalized Verma module in the above proposition is called an  elementary homomorphism (\cite{[M2]}, \cite{[M3]}).

The goal of this article is to prove the following result.

\begin{thm}
Non-zero homomorphisms between scalar generalized Verma modules of $\gl(n, \cpx)$ are compositions of elementary homomorphisms.
\end{thm}

The above theorem confirms Conjecture A in \cite{[M3]} for ${\mathfrak gl}(n, \cpx)$.

As we explained in \cite{[M3]} 2.2, Soergel's result (\cite{[So]} Theorem 11 also see \cite{[Hum]} 13.13) implies that the above theorem is reduced to the following integral weight case.

\begin{thm}
Let  $\mu_1,....\mu_k\in\itg$ and  $\nu_1,....\nu_k\in\itg$.  Then, any non-zero homomorphism of  $ M_\T([\nu_1,...\nu_k])$ to $M_\T([\mu_1,...\mu_k])$ is a composition of elementary homomorphisms.
\end{thm}

\remark  In fact, the argument of our proof also works for non-integral weights.  So, the above-mentioned reduction is not essential for us.
However, without the reduction, the proof requires more complicated notations.

\subsection{Reformulation in terms of a Bruhat ordering}

First, we show the following result.

\begin{lem}
Let $\mu,\nu\in{}^\circ\PS$ and we assume that $M_\T(\nu)\subseteq M_\T(\mu)$.  Then, there exists some $x\in W(\T)$ such that $x\nu=\mu$.
\end{lem}

\proof
We put $I_1=\Ann_{U(\gggg)}(M_{\T}(\mu))$ and $I_2=\Ann_{U(\gggg)}(M_{\T}(\nu))$.
$I_1$ and $I_2$ are primitive ideals. (See \cite{[J]} 15.6 Korollar.)
We see that $M_{\T}(\nu)\subseteq M_{\T}(\mu)$ implies
$I_1\subseteq I_2$.
Since the Gelfand-Kirillov dimension of  $M_{\T}(\mu)$ and $M_{\T}(\mu)$ are same,
we have $I_1=I_2$ from \cite{[BK]} 3.6.Korollar.
Then, we obtain the the lemma from \cite{[BJ]} 5.10 Korollar.
\,\,\,\, $\Box$

Let $\lambda_1,...,\lambda_k\in\itg$.  We call $[\lambda_1,...,\lambda_k]$ $\T$-antidominant, if $\lambda_p\leqslant\lambda_q$ for all $1\leqslant p<q\leqslant k$ such that $n_p=n_q$. 

We easily see the following result.
\begin{lem}
For any $\mu_1,...,\mu_k\in\itg$, there is a unique $\T$-antidominant $[\lambda_1,...,\lambda_k]$ which is $W(\T)$-conjugate to $[\mu_1,...,\mu_k]$.
\end{lem}

We fix a  $\T$-antidominant $\lambda=[\lambda_1,...,\lambda_k]\in {}^\circ\PS$.
We introduce a Bruhat ordering on $W(\T)\lambda=\{w\lambda\mid w\in W(\T)\}$ as follows.
First, for $\mu=[\mu_1,...\mu_k],\nu=[\nu_1,...,\nu_k]\in W(\T)\lambda$, we write $\nu\uparrow\mu$, if there exist $1\leqslant p<q\leqslant k$ such that $n_p=n_q$, $\nu=\sigma_{p,q}\mu$, and $\lambda_p-\lambda_q\in\nat$.
For $\mu,\nu\in W(\T)\lambda$, we write $\nu\leqslant_\T \mu$, if there exists some finite sequence $\eta_1,...,\eta_h \in W(\T)\lambda$ such that $\nu\uparrow\eta_1\uparrow\eta_2\uparrow\cdots\uparrow\eta_h\uparrow\mu$.

As in the case of Verma modules (cf.\ \cite{[BGG]}), we can rephrase Theorem 2.3.3 as follows.
\begin{thm}
Let  $\lambda_1,...,\lambda_k\in\itg$ be such that $\lambda=[\lambda_1,...,\lambda_k]$ is $\T$-antidominant.
Let $\mu,\nu\in W(\T)\lambda$.  Then the following (a),(b) are equivalent.

(a) \,\,\, $M_\T(\nu)\subseteq M_\T(\mu)$.

(b) \,\,\, $\nu\leqslant_\T \mu$.
\end{thm}

We immediately see that (b) implies (a) from Proposition 2.3.1.
So, we have only to show that (a) implies (b).

\setcounter{section}{3}
\setcounter{subsection}{0}

\section*{\S\,\,3.\,\,\,\,Proof of the main result }
\subsection{Tableau description of the Bruhat orderings}

We fix  $\lambda_1,...,\lambda_k\in\itg$  such that $\lambda=[\lambda_1,...,\lambda_k]$ is $\T$-antidominant.

For $\mu=[\mu_1,...,\mu_k]\in W(T)\lambda$ and $r\in\Upsilon$ , we write 
$$\mu_{(r)}=\sum_{j\in I_r}\sum_{s=1}^{n_j}(\lambda_j-s+1)e_{n_{j-1}^\ast+s}.$$
Then, we immediately see $\mu=\sum_{r\in\Upsilon}\mu_{(r)}$.

Using the identification $(\sharp)$ in 2.2, we write an element of $W(\T)$ as $(w_r)_{r\in \Upsilon}$, where $w_r\in W(\T)_r$.
Then, we easily see $w\mu=\sum_{r\in\Upsilon}w_r\mu_{(r)}$ for all $\mu\in W(T)\lambda$.
Hence, $w\lambda\rightsquigarrow (w_r\lambda_{(r)})_{r\in\Upsilon}$ gives
\begin{align*}
W(\T)\cong \prod_{r\in\Upsilon}W(\T)_r\lambda_{(r)}. \tag{$\sharp\sharp$}
\end{align*}
For $\mu_{(r)}, \nu_{(r)}\in W(\T)_r\lambda_{(r)}$, we write $\nu_{(r)}\uparrow_r\mu_{(r)}$ if there exists some $p,q\in I_r$ such that $p<q$, $\sigma_{p,q}$, and $\mu_p-\mu_q\in\nat$.
The Bruhat ordering $\leqslant_r$ on $W(\T)_r\lambda_{(r)}$ is the partial order generated by $\uparrow_r$. 
So, we see that the following ($\triangle$) holds.
\begin{align*}
\hbox{$\nu\leqslant_\T\mu$ if and only if $\nu_{(r)}\leqslant_r\mu_{(r)}$ for all $r\in\Upsilon$.}\tag{$\triangle$}
\end{align*}

Next, we explain a tableau description of $(W(\T)_r\lambda_{(r)},\leqslant_r)$ in \cite{[P]}.
(It seems that the idea of the description goes back to \cite{[E]}.)
We fix $r\in\Upsilon$ and  $\lambda_1,...,\lambda_k\in\itg$ be such that $\lambda=[\lambda_1,...,\lambda_k]$ be $\T$-antidominant.
We put $M_r=\card\{\lambda_i\mid i\in I_r\}$.
Let $\eta_1,...,\eta_{M_r}\in \itg$ be such that $\{\eta_1<\eta_2<\cdots<\eta_{M_r}\}=\{\lambda_i\mid i\in I_r\}$.

For $\mu\in W(\T)\lambda$ and $1\leqslant i\leqslant M_r$, we put 
$$T_{(r)}^i[\mu]=\{j\in I_r\mid \mu_j\geqslant \eta_i\}.$$
We put $N_r^i=\card T_{(r)}^i[\lambda]$ for $1\leqslant i\leqslant M_r$.
We immediately see $\card T_{(r)}^i[\mu]=N_r^i$.
We define $t_{(r),j}^i[\mu]$ for  $1\leqslant i\leqslant M_r$ and $1\leqslant j\leqslant N_r^i$ as follows.
$$\left\{t_{(r),1}^i[\mu]<\cdots<t_{(r),N_r^i}^i[\mu]\right\}=T_{(r)}^i[\mu].$$

\begin{lem} \, (\cite{[P]} Theorem 5A)

We fix  $\lambda_1,...,\lambda_k\in\itg$  such that $\lambda=[\lambda_1,...,\lambda_k]$ is $\T$-antidominant.
For $\mu,\nu\in W(\T)\lambda$ and $r\in \Upsilon$, the following (1) and (2) are equivalent.
\begin{itemize}
\item[(1)] \,\, $\nu_{(r)}\leqslant_r\mu_{(r)}$.
\item[(2)] \,\, $t_{(r),j}^i[\nu]\geqslant t_{(r),j}^i[\mu] $ for all $1\leqslant i\leqslant M_r$ and $1\leqslant j\leqslant N_r^i$.
\end{itemize}
 
\end{lem}

\remark  We remark that the positive roots in \cite{[P]} are the negative roots in our setting.  This affects the statement of the above lemma.

Using the tableau description, we have the following result.
\begin{lem}
 Let $\lambda_1,...\lambda_k\in\itg$ be such that $\lambda=[\lambda_1,...,\lambda_k]$ is $\T$-antidominant. Let $x,y\in W(\T)$ and we write  $\nu=[\nu_1,...,\nu_k]=x\lambda, \mu=[\mu_1,...,\mu_k]=y\lambda$. 
We assume that $\nu\not\leqslant_\T\mu$.
Then there exists some $c\in\itg$ and $1\leqslant s\leqslant k$ satisfying the following (1) and (2).
\begin{itemize}
\item[(1)] \,\, $\nu_s\geqslant c>\mu_s$.
\item[(2)] \,$\card\{i\mid 1\leqslant i<s, n_i=n_s, \nu_i\geqslant c\}=\card\{i\mid  1\leqslant i<s, n_i=n_s, \mu_i\geqslant c\}.$
\end{itemize}
\end{lem}

\proof
From the above definitions, we see that there exists some $r\in I_r$ such that $\nu_{(r)}\leqslant_r\mu_{(r)}$.
From Lemma 3.1.1, there exist some $1\leqslant i\leqslant M_r$ and $1\leqslant j\leqslant N_r^i$ such that  $t_{(r),j}^i[\nu]< t_{(r),j}^i[\mu]$.  We choose $1\leqslant a\leqslant M_r$ and $1\leqslant b\leqslant M_r^i$ such that the following (A),(B), and (C) hold.
\begin{itemize}
\item[(A)] \,\, $t_{(r),b}^a[\nu]< t_{(r),b}^v[\mu]$.
\item[(B)] \,\, $t_{(r),j}^{i}[\nu]\geqslant t_{(r),j}^{i}[\mu]$ for all $a<i\leqslant M_r$ and  $1\leqslant j\leqslant N_r^{a}$.
\item[(C)] \,\, $t_{(r),j}^{a}[\nu]\geqslant t_{(r),j}^{a}[\mu]$ for all  $1\leqslant j<b$.
\end{itemize}
We put $s= t_{(r),b}^a[\nu]$ and $c=\eta_a$.  From the definition, we have $\nu_s\geqslant c$.
From the definition, we also see:
$$\{i\mid 1\leqslant i<s, n_i=n_s, \nu_i\geqslant c\}=\{t_{(r),j}^{a}[\nu]\mid 1\leqslant j<b\}.$$
From (A) and (C) above, we have
$$\{i\mid 1\leqslant i<s, n_i=n_s, \mu_i\geqslant c\}=\{t_{(r),j}^{a}[\mu]\mid 1\leqslant j<b\}.$$ 
Hence, we have (2).

From (A) and (C), we also have
$$ t_{(r),b-1}^{a}[\mu]\leqslant t_{(r),b-1}^{a}[\nu]<t_{(r),b}^{a}[\nu]=s<t_{(r),b}^{a}[\mu].$$
Since there is no $j\in T_{(r)}^{a}[\mu]$ such that $t_{(r),b-1}^{a}[\mu]<j<t_{(r),b}^{a}[\mu]$, we have $s\not\in T_{(r)}^{a}[\mu]$.
This means that $\mu_s<\eta_a=c$.  So, we have (1). \,\,\,\,$\Box$

\begin{lem}
 Let $\lambda_1,...\lambda_k\in\itg$ be such that $\lambda=[\lambda_1,...,\lambda_k]$ is $\T$-antidominant. Let $x,y\in W(\T)$ and we write  $\nu=[\nu_1,...,\nu_k]=x\lambda, \mu=[\mu_1,...,\mu_k]=y\lambda \in W(\T)$. 
We assume that $\nu\not\leqslant_\T\mu$.
We choose $c\in\itg$ and $1\leqslant s\leqslant k$ as in Lemma 3.1.2 and define $\bar{\mu_1},...,\bar{\mu}_k$ and $\bar{\nu_1},...,\bar{\nu}_k$ as follows.
$$\bar{\mu}_i=
\begin{cases}
c  & \hbox{if $n_i>n_s$}\\
c-1  & \hbox{if $n_i<n_s$}\\
c  & \hbox{if $n_i=n_s$ and $\mu_i\geqslant c$}\\
c-1  & \hbox{if $n_i=n_s$ and $\mu_i< c$}
\end{cases}.
$$

$$\bar{\nu}_i=
\begin{cases}
c  & \hbox{if $n_i>n_s$}\\
c-1  & \hbox{if $n_i<n_s$}\\
c  & \hbox{if $n_i=n_s$ and $\nu_i\geqslant c$}\\
c-1  & \hbox{if $n_i=n_s$ and $\nu_i< c$}
\end{cases}.
$$
Then we have
\begin{itemize}
\item[(1)] \,\,  $[\bar{\nu}_1,...,\bar{\nu}_k]=xy^{-1}[\bar{\mu}_1,...,\bar{\mu}_k]$.
\item[(2)] \,\,  $M_\T([\bar{\nu}_1,...,\bar{\nu}_k])\not\subseteq  M_\T([\bar{\mu}_1,...,\bar{\mu}_k]).$
\end{itemize}
\end{lem}

\proof
$ [\bar{\nu}_1,...,\bar{\nu}_k]=xy^{-1}[\bar{\mu}_1,...,\bar{\mu}_k]$ easily follows from $[{\nu}_1,...,{\nu}_k]=xy^{-1}[{\mu}_1,...,{\mu}_k]$ and the definition.

We prove (2).
We assume $ M_\T([\bar{\nu}_1,...,\bar{\nu}_k])\subseteq  M_\T([\bar{\mu}_1,...,\bar{\mu}_k])$ and deduce a contradiction.
We easily see that there are some non-negative integers $\ell_1,...\ell_n$ such that $[\bar{\mu}_1,...,\bar{\mu}_k]-[\bar{\nu}_1,...,\bar{\nu}_k]=\sum_{i=1}^{n-1}\ell_i(e_i-e_{i+1})$.  Since we have 
$$[\bar{\mu}_1,...,\bar{\mu}_k]-[\bar{\nu}_1,...,\bar{\nu}_k]=\sum_{i=1}^k\left((\bar{\mu}_i-\bar{\nu}_i)\sum_{j=1}^{n_i}e_{n^\ast_{i-1}+j}\right),$$
we see 
$$\sum_{i=1}^k\left((\bar{\mu}_i-\bar{\nu}_i)\sum_{j=1}^{n_i}e_{n^\ast_{i-1}+j}\right)=\sum_{i=1}^{n}(\ell_i-\ell_{i-1})e_i.$$
Here, we put $\ell_0=\ell_n=0$.
From Lemma 3.1.2, we have 
$$\ell_{n^\ast_{s-1}}=\sum_{i=1}^{n^\ast_{s-1}}(\ell_i-\ell_{i-1})=\sum_{i=1}^{s-1}n_i(\bar{\mu}_i-\bar{\nu}_i)=0.$$
Hence, we have
$$\ell_{n^\ast_{s-1}+1}=\ell_{n^\ast_{s-1}+1}-\ell_{n^\ast_{s-1}}=\bar{\mu}_{s}-\bar{\nu}_{s}=(c-1)-c=-1.$$
It contradicts the non-negativity of $\ell_{n^\ast_{s-1}+1}$.  \,\,\,\, $\Box$

\subsection{Translations in a mediocre region}

The material in this subsection is found more or less in \cite{[K]}, \cite{[T]}, \cite{[VU]}, \cite{[VI]}, \cite{[Vd]}.
(Our usage of ``a mediocre region" is not necessarily ``the mediocre range for $\pps$" in \cite{[T]} Definition 3.4, but it is the mediocre 
 range for some parabolic subalgebra with a Levi part $\lls$.)

For $1\leqslant i\leqslant k$, we put $f_i=\sum_{j=1}^{n_i}e_{n^\ast_{i-1}+j}$.
For $\lambda_1,...,\lambda_k\in\cpx$, we easily see 
$$[\lambda_1,...,\lambda_k]\pm f_i=[\lambda_1,...,\lambda_{i-1},\lambda_i\pm 1,\lambda_{i+1},...,\lambda_k].$$
For $1\leqslant i\leqslant k$ and $1\leqslant r\leqslant n_i$, we also put as follows.
$$\overline{f}_i(r)=\sum_{j=1}^re_{n^\ast_{i-1}+j},$$
$$\underline{f}_i(r)=\sum_{j=1}^re_{n^\ast_{i}-j+1}.$$
So,we see $f_i=\overline{f}_i(n_i)=\underline{f}_i(n_i).$
We also put $\overline{f}_i(0)=\underline{f}_i(0)=0$.

We denote by $\square_m$ (resp.\ $\square^\ast_m$) the natural representation(resp.\ the dual of natural representation)  of ${\mathfrak gl}(m,\cpx)$ and denote by $\wedge^r\square_m$ (resp.\  $\wedge^r\square_m^\ast$) its $r$-th wedge product representation.  Since $\llll_\T\cong \gl(n_1,\cpx)\oplus\cdots\oplus\gl(n_k,\cpx)$, we may regard an external tensor product $\wedge^{r_1}\square_{n_1}\boxtimes\cdots\boxtimes\wedge^{r_k}\square_{n_k}$ as an $\llll_\T$-module.
The following result is well-known.
\begin{lem}
For $1\leqslant r\leqslant n$, we have
$$\wedge^r\square_n|_{\llll_\T}\cong\bigoplus_{r_1+\cdots+r_k=r}\wedge^{r_1}\square_{n_1}\boxtimes\cdots\boxtimes\wedge^{r_k}\square_{n_k},$$
$$\wedge^r\square_n^\ast|_{\llll_\T}\cong\bigoplus_{r_1+\cdots+r_k=r}\wedge^{r_1}\square_{n_1}^\ast\boxtimes\cdots\boxtimes\wedge^{r_k}\square_{n_k}^\ast.$$
\end{lem}
Hence, we easily see the following result.
\begin{lem}
For $1\leqslant r\leqslant n$, $\wedge^r\square_n|_{\ppp_\T}$ has a filtration of $\pps$-submodules such that the set of its successive quotients is 
$$\{\wedge^{r_1}\square_{n_1}\boxtimes\cdots\boxtimes\wedge^{r_k}\square_{n_k}\mid r_1+\cdots+r_k=r, 0\leqslant r_i\leqslant n_i \,\,\, (1\leqslant i\leqslant k) \}.$$
Here, we regard $\wedge^{r_1}\square_{n_1}\boxtimes\cdots\boxtimes\wedge^{r_k}\square_{n_k}$ as a $\pps$-module on which $\nns$ acts trivially.
Similarly, 
 $\wedge^r\square^\ast_n|_{\ppp_\T}$ has a filtration of $\pps$-submodules such that the set of its successive quotients is 
$$\{\wedge^{r_1}\square_{n_1}^\ast\boxtimes\cdots\boxtimes\wedge^{r_k}\square_{n_k}^\ast\mid r_1+\cdots+r_k=r, 0\leqslant r_i\leqslant n_i \,\,\, (1\leqslant i\leqslant k)\}.$$
\end{lem}
From an infinitesimal version of Mackey tensor product theorem, for $\lambda\in\PS$, we have
$$ M_\T(\lambda)\otimes\wedge^r\square_n\cong U(\ggg)\otimes_{U(\pps)}(E_\T(\lambda)\otimes\wedge^r\square_n|_{\pps}).$$
Hence, we have the following result.
\begin{lem}
Let $\lambda\in{}^\circ\PS$.
Then, $M_\T(\lambda)\otimes\wedge^r\square_n$ has a filtration of $\ggg$-submodules such that the set of its successive quotients is 
$$\{M_\T(\lambda+\overline{f}_1(r_1)+\cdots+\overline{f}_k(r_k))\mid  r_1+\cdots+r_k=r, 0\leqslant r_i\leqslant n_i \,\,\, (1\leqslant i\leqslant k)\}.$$
Similarly,  $M_\T(\lambda)\otimes\wedge^r\square_n^\ast$ has a filtration of $\ggg$-submodules such that the set of its successive quotients is 
$$\{M_\T(\lambda-\underline{f}_1(r_1)-\cdots-\underline{f}_k(r_k))\mid  r_1+\cdots+r_k=r, 0\leqslant r_i\leqslant n_i \,\,\, (1\leqslant i\leqslant k)\}.$$

\end{lem}
For $1\leqslant i_1,i_2,...,i_s\leqslant k$ and $\lambda\in{}^\circ\PS$, we define translation functors $T^{\lambda\pm \sum_{r=1}^s f_{i_r}}_\lambda: \mca_\lambda\rightarrow\mca_{\lambda\pm \sum_{r=1}^s f_{i_r}}$ as follows.
$$T^{\lambda+ \sum_{r=1}^s f_{i_r}}_\lambda(M)={\bf P}_{\lambda+ \sum_{r=1}^s f_{i_r}}(M\otimes\wedge^{\sum_{r=1}^s n_{i_r}}\square_n),$$
$$T^{\lambda- \sum_{r=1}^s f_{i_r}}_\lambda(M)={\bf P}_{\lambda- \sum_{r=1}^s f_{i_r}}(M\otimes\wedge^{\sum_{r=1}^s n_{i_r}}\square^\ast_n).$$

The following result is easy.

\begin{lem}
Let $\lambda\in{}^\circ\PS$ and $w\in W(\T)$.  Then, we have  $T^{\lambda\pm \sum_{r=1}^s f_{i_r}}_\lambda=T^{w\lambda\pm \sum_{r=1}^s f_{\bar{w}(i_r)}}_{w\lambda}$.
\end{lem}

For $\lambda\in\hhd$ and $g\in\itg$, we put $\|\lambda; g\|=\card\{i\mid 1\leqslant i\leqslant n, \langle\lambda,e_i\rangle=g\}$.
Obviously, for $w\in W$, we have
\begin{align*}\| w\lambda;g\|=\|\lambda;g\|. \tag{$\dagger$}\end{align*}

For $\lambda_1,...,\lambda_k\in\itg$, we put
$$\Phi([\lambda_1,...,\lambda_k];g)=\{i\mid 1\leqslant i\leqslant k, \lambda_i\geqslant g\geqslant\lambda_i-n_i+1\}.$$
From the definition, we easily see 
$$\| [\lambda_1,...,\lambda_k];g\|=\card\Phi([\lambda_1,...,\lambda_k];g).$$
We also put
$$\overline{\Psi}([\lambda_1,...,\lambda_k];g)=\{i\mid 1\leqslant i\leqslant k, \lambda_i=g\},$$
$$\underline{\Psi}([\lambda_1,...,\lambda_k];g)=\{i\mid 1\leqslant i\leqslant k, \lambda_i-n_i+1=g\},$$

We have the following result.
\begin{lem}
Let $\lambda_1,...,\lambda_k\in\itg$ and let  $g\in\itg$. We put $\lambda=[\lambda_1,...,\lambda_k]$.
We fix $S\subseteq \overline{\Psi}(\lambda;g)$.
We assume that the following condition (a) holds.
\begin{itemize}
\item[(a)]  $\min\{\lambda_i-n_{i}+1\mid i\in S\}>\max\{\lambda_j-n_j+1\mid j\in \Phi(\lambda;g)-S\}$.
\end{itemize}
Then, we have
$$T_\lambda^{\lambda-\sum_{j\in S}f_j}(M_\T(\lambda))\cong M_\T\left(\lambda-\sum_{j\in S}f_j\right).$$
Here, if $\Phi(\lambda;g)-S=\emptyset$, we regard (a) as an empty condition.
\end{lem}
\proof
We put $r=\sum_{j\in S}n_j.$
From Lemma 3.2.3, we see that $T_\lambda^{\lambda-\sum_{j\in S}f_j}(M_\T(\lambda))$ has a filtration of $\ggg$-submodules such that the set of its successive quotients is
\begin{multline*}
{\bf SQ}=\left\{M_\T\left(\lambda-\sum_{i=1}^k\underline{f}_i(r_i)\right)\right| \sum_{i=1}^k r_i=r, 0\leqslant r_i\leqslant n_i \, (1\leqslant i\leqslant k),  \\ \left.\exists w\in W \left[\lambda-\sum_{i=1}^k\underline{f}_i(r_i)=w\left(\lambda-\sum_{j\in S}f_j\right)\right]\right\}.
\end{multline*}
So, clearly, we have $M_\T\left(\lambda-\sum_{j\in S}f_j\right)\in{\bf SQ}$.

On the other hand, we assume that $M_\T\left(\lambda-\sum_{i=1}^k\underline{f}_i(r_i)\right)\in {\bf SQ}$.
Put $\ell=\card S.$ From $(\dagger)$, we have
$$\ell=\|\lambda;g\|-\left\|\lambda-\sum_{j\in S}f_j;g\right\|=\|\lambda;g\|-\left\|\lambda-\sum_{i=1}^k\underline{f}_i(r_i);g\right\|.$$
Hence, there exist integers $1\leqslant h_1<\cdots<h_\ell\leqslant k$ such that $r_{h_i}=g-\lambda_{h_i}+n_{h_i}$ \,\,\, $(1\leqslant i\leqslant\ell)$.
If $\lambda-\sum_{i=1}^k\underline{f}_i(r_i)\neq\lambda-\sum_{j\in S}f_j$, there exists some $1\leqslant s\leqslant k$ such that $h_s\in \Phi(\lambda;g)-S$.
Since $g=\lambda_i$ for $i\in S$, the condition (a) can be rephrased as follows.
$$\max\{n_{i}\mid i\in S\}<\min\{g-\lambda_j+n_j\mid j\in \Phi(\lambda;g)-S \}$$
Hence, we have
$$r=\sum_{j\in S}n_j<\sum_{i=1}^\ell r_{h_i}\leqslant\sum_{j=1}^k r_j=r.$$
Hence, we obtain a contradiction. \,\,\,\,$\Box$

Similarly, we also have the following lemma.
\begin{lem}
Let $\lambda_1,...,\lambda_k\in\itg$ and let  $g\in\itg$. We put $\lambda=[\lambda_1,...,\lambda_k]$.
We fix $S\subseteq \underline{\Psi}(\lambda;g)$.
We assume that the following condition (b) holds.
\begin{itemize}
\item[(b)]  $\max\{\lambda_i\mid i\in S\}<\min\{\lambda_j\mid j\in \Phi(\lambda;g)-S \}$.
\end{itemize}
Then, we have
$$T_\lambda^{\lambda+\sum_{j\in S}f_j}(M_\T(\lambda))\cong M_\T\left(\lambda+\sum_{j\in S}f_j\right).$$
Here, if $\Phi(\lambda;g)-S=\emptyset$, we regard (b) as an empty condition.
\end{lem}
\subsection{Proof of Theorem 2.4.3}
As we remarked in the last paragraph of 2.4, we have only to show that (a) in 2.4.3 implies (b) in 2.4.3.
So, we  let  $\lambda_1,...,\lambda_k\in\itg$ be such that $[\lambda_1,...,\lambda_k]$ is $\T$-antidominant.
Let $x,y\in W(\T)$. 
Put $\nu=[\nu_1,...,\nu_k]=x[\lambda_1,...,\lambda_k]$ and $\mu=[\mu_1,...,\mu_k]=y[\lambda_1,...,\lambda_k]$, and $w=xy^{-1}\in W(\T)$. 
We denote by $\bar{w}$ the element in ${\mathfrak S}_k$ corresponding to $w$  via ($\flat$) in 2.2.

 We assume that  $M_\T([\nu_1,...,\nu_k])\subseteq M_\T([\mu_1,...,\mu_k])$ and
 $\nu\not\leqslant_\T \mu$ and deduce a contradiction.

From Lemma 3.1.2, there exists some $c\in\itg$ and $1\leqslant s\leqslant k$ satisfying the conditions (1) and (2) in the statement of Lemma 3.1.2.
We define $[\bar{\mu}_1,...,\bar{\mu}_k]$ and  $[\bar{\nu}_1,...,\bar{\nu}_k]$ as in the statement of Lemma 3.1.3.
Since the translation functors are exact, they map injective homomorphisms to injective homomorphisms.
The idea of the proof is, applying translation functors several times, to derive ``$M_\T([\bar{\nu}_1,...,\bar{\nu}_k])\subseteq M_\T([\bar{\mu}_1,...,\bar{\mu}_k])$" from ``$M_\T([\nu_1,...,\nu_k])\subseteq M_\T([\mu_1,...,\mu_k])$".  We divide this procedure into several steps.

\medskip

{\it Step 1}  \,\,\,\, For $d\geqslant 0$, we define $\widehat{\mu}^{(d)}=[\widehat{\mu}_1^{(d)},...,\widehat{\mu}_k^{(d)}]$ and  $\widehat{\nu}^{(d)}=[\widehat{\nu}_1^{(d)},...,\widehat{\nu}_k^{(d)}]$ as follows.

$$\widehat{\mu}_i^{(d)}=
\begin{cases}
c+d  & \hbox{if $\mu_i\geqslant c+d$}\\
\mu_i  & \hbox{if $\mu_i< c+d$}
\end{cases}.
$$

$$\widehat{\nu}_i^{(d)}=
\begin{cases}
c+d  & \hbox{if $\nu_i\geqslant c+d$}\\
\nu_i  & \hbox{if $\nu_i< c+d$}
\end{cases}.
$$

We easily see $\widehat{\nu}^{(d)}=w\widehat{\mu}^{(d)}$.
Obviously, we have $\overline{\Psi}(\widehat{\nu}^{(d)};c+d)=\Phi(\widehat{\nu}^{(d)};c+d)$ and   $\overline{\Psi}(\widehat{\mu}^{(d)};c+d)=\Phi(\widehat{\mu}^{(d)};c+d)$.

Hence, for all $d>0$, we have
$$\widehat{\mu}^{(d)}-\sum_{j\in\overline{\Psi}(\widehat{\mu}^{(d)};c+d)}f_j=\widehat{\mu}^{(d-1)}.$$
Therefore, from Lemma 3.2.5, we have
$$T_{\widehat{\mu}^{(d)}}^{\widehat{\mu}^{(d-1)}}\left(M_\T\left(\widehat{\mu}^{(d)}\right)\right)\cong M_\T\left(\widehat{\mu}^{(d-1)}\right).$$
Similarly, we have
$$T_{\widehat{\nu}^{(d)}}^{\widehat{\nu}^{(d-1)}}\left(M_\T\left(\widehat{\nu}^{(d)}\right)\right)\cong M_\T\left(\widehat{\nu}^{(d-1)}\right).$$
From Lemma 3.2.4, we have $$T_{\widehat{\mu}^{(d)}}^{\widehat{\mu}^{(d-1)}}=T_{w\widehat{\mu}^{(d)}}^{w\widehat{\mu}^{(d-1)}}=T_{\widehat{\nu}^{(d)}}^{\widehat{\nu}^{(d-1)}}.$$
If $d$ is sufficiently large, clearly we have $\widehat{\mu}^{(d)}=\mu$ and   $\widehat{\nu}^{(d)}=\nu$.
Since translation functors are exact, we have
 $$M_\T\left(\widehat{\nu}^{(0)}\right)\subseteq
 M_\T\left(\widehat{\mu}^{(0)}\right)$$
 from $M_\T\left({\nu}\right)\subseteq
 M_\T\left({\mu}\right).$ 

\medskip

{\it Step 2}  \,\,\,\, We define $\widehat{\mu}=[\widehat{\mu}_1,...,\widehat{\mu}_k]$ and  $\widehat{\nu}=[\widehat{\nu}_1,...,\widehat{\nu}_k]$ as follows.

$$\widehat{\mu}_i=
\begin{cases}
c  & \hbox{if $\mu_i\geqslant c$ and $n_i\geqslant n_s$}\\
c-1  & \hbox{if $\mu_i\geqslant c$ and $n_i< n_s$}\\
\mu_i  & \hbox{if $\mu_i< c$}
\end{cases},
$$
$$\widehat{\nu}_i=
\begin{cases}
c  & \hbox{if $\nu_i\geqslant c$ and $n_i\geqslant n_s$}\\
c-1  & \hbox{if $\nu_i\geqslant c$ and $n_i< n_s$}\\
\nu_i  & \hbox{if $\nu_i< c$}
\end{cases}.
$$
We put $S_1=\{i\mid 1\leqslant i\leqslant k, \mu_i\geqslant c, n_i<n_s\}.$
Then, we easily see $S_1=\{i\in \overline{\Psi}(\widehat{\mu}^{(0)};c)\mid n_i<n_s\}$ and $\widehat{\mu}=\widehat{\mu}^{(0)}-\sum_{i\in S_1}f_i$.
From  Lemma 3.2.5, we have
$$T_{\widehat{\mu}^{(0)}}^{\widehat{\mu}}\left(M_\T\left(\widehat{\mu}^{(0)}\right)\right)\cong M_\T\left(\widehat{\mu}\right).$$
Similarly, we have
$$T_{\widehat{\nu}^{(0)}}^{\widehat{\nu}}\left(M_\T\left(\widehat{\nu}^{(0)}\right)\right)\cong M_\T\left(\widehat{\nu}\right).$$
Since $w\widehat{\mu}=\widehat{\nu}$, Lemma 3.2.4 implies
$$T_{\widehat{\mu}^{(0)}}^{\widehat{\mu}}=T_{w\widehat{\mu}^{(0)}}^{w\widehat{\mu}}=T_{\widehat{\nu}^{(0)}}^{\widehat{\nu}}.$$
Since translation functors are exact, we have 
$$M_\T\left(\widehat{\nu}\right)\subseteq
 M_\T\left(\widehat{\mu}\right)$$
 from $M_\T\left(\widehat{\nu}^{(0)}\right)\subseteq
 M_\T\left(\widehat{\mu}^{(0)}\right).$ 

\medskip

{\it Step 3}  \,\,\,\, We define an integer $b$ as follows.
$$b=
\begin{cases}
c-\max\{n_i\mid 1\leqslant i\leqslant k\}+1 & \hbox{if there is some $1\leqslant i\leqslant k$ such that $n_i>n_s$,}\\
c-n_s  & \hbox{otherwise.}
\end{cases}.
$$
For $d\geqslant 0$,  we define $\widetilde{\mu}^{(d)}=[\widetilde{\mu}^{(d)}_1,...,\widetilde{\mu}^{(d)}_k]$ and  $\widetilde{\nu}^{(d)}=[\widetilde{\nu}^{(d)}_1,...,\widetilde{\nu}^{(d)}_k]$ as follows.

$$\widetilde{\mu}^{(d)}_i=
\begin{cases}
c  & \hbox{if $\mu_i\geqslant c$ and $n_i\geqslant n_s$}\\
c-1  & \hbox{if $\mu_i\geqslant c$ and $n_i< n_s$}\\
\mu_i  & \hbox{if $b-d+n_i-1\leqslant\mu_i< c$}\\
b-d+n_i-1  & \hbox{if $\mu_i<b-d+n_i-1$,}
\end{cases}
$$
$$\widetilde{\nu}^{(d)}_i=
\begin{cases}
c  & \hbox{if $\nu_i\geqslant c$ and $n_i\geqslant n_s$}\\
c-1  & \hbox{if $\nu_i\geqslant c$ and $n_i< n_s$}\\
\nu_i  & \hbox{if $b-d+n_i-1\leqslant\nu_i< c$}\\
b-d+n_i-1  & \hbox{if $\nu_i<b-d+n_i-1$.}
\end{cases}
$$
Obviously, we have $\widetilde{\nu}^{(d)}=w\widetilde{\mu}^{(d)}$, $\underline{\Psi}(\widetilde{\nu}^{(d)};b-d)=\Phi(\widetilde{\nu}^{(d)};b-d)$, and   $\underline{\Psi}(\widetilde{\mu}^{(d)};b-d)=\Phi(\widetilde{\mu}^{(d)};b-d)$.

Hence, for all $d>0$, we have
$$\widetilde{\mu}^{(d)}+\sum_{j\in\underline{\Psi}(\widetilde{\mu}^{(d)};b-d)}f_j=\widetilde{\mu}^{(d-1)}.$$
Therefore, from Lemma 3.2.6, we have
$$T_{\widetilde{\mu}^{(d)}}^{\widetilde{\mu}^{(d-1)}}\left(M_\T\left(\widetilde{\mu}^{(d)}\right)\right)\cong M_\T\left(\widetilde{\mu}^{(d-1)}\right).$$
Similarly, we have
$$T_{\widetilde{\nu}^{(d)}}^{\widetilde{\nu}^{(d-1)}}\left(M_\T\left(\widetilde{\nu}^{(d)}\right)\right)\cong M_\T\left(\widetilde{\nu}^{(d-1)}\right).$$

From Lemma 3.2.4, we have $$T_{\widetilde{\mu}^{(d)}}^{\widetilde{\mu}^{(d-1)}}=T_{w\widetilde{\mu}^{(d)}}^{w\widetilde{\mu}^{(d-1)}}=T_{\widetilde{\nu}^{(d)}}^{\widetilde{\nu}^{(d-1)}}.$$
If $d$ is sufficiently large, clearly we have $\widetilde{\mu}^{(d)}=\widehat{\mu}$ and   $\widetilde{\nu}^{(d)}=\widehat{\nu}$.
Since translation functors are exact, we have
 $$M_\T\left(\widetilde{\nu}^{(0)}\right)\subseteq
 M_\T\left(\widetilde{\mu}^{(0)}\right)$$
 from $M_\T\left(\widehat{\nu}\right)\subseteq
 M_\T\left(\widehat{\mu}\right).$ 

\medskip

{\it Step 4}  \,\,\,\,
We skip this step if there is no $1\leqslant i\leqslant k$ such that $n_s<n_i$.
So, we assume $n_s<\max\{n_i\mid 1\leqslant i\leqslant k\}$.
Let $b$ be as in Step 3. Namely, $b=c-\max\{n_i\mid 1\leqslant i\leqslant k\}+1$.
For $0\leqslant d\leqslant\max\{n_i\mid 1\leqslant i\leqslant k\}-n_s-1=c-n_s-b$, we define   $\breve{\mu}^{(d)}=[\breve{\mu}^{(d)}_1,...,\breve{\mu}^{(d)}_k]$ and  $\breve{\nu}^{(d)}=[\breve{\nu}^{(d)}_1,...,\breve{\nu}^{(d)}_k]$ as follows.

$$\breve{\mu}^{(d)}_i=
\begin{cases}
c  & \hbox{if $\mu_i\geqslant c$ and $n_i\geqslant n_s$}\\
c-1  & \hbox{if $\mu_i\geqslant c$ and $n_i< n_s$}\\
\mu_i  & \hbox{if $c-n_s+n_i-1-d\leqslant\mu_i< c$}\\
c  & \hbox{if $\mu_i<c-n_s+n_i-1-d$ and $n_i\geqslant n_s+d+1$}\\
c-n_s+n_i-1-d  & \hbox{if $\mu_i<c-n_s+n_i-1-d$ and $n_i< n_s+d+1$,}
\end{cases}
$$
$$\breve{\nu}^{(d)}_i=
\begin{cases}
c  & \hbox{if $\nu_i\geqslant c$ and $n_i\geqslant n_s$}\\
c-1  & \hbox{if $\nu_i\geqslant c$ and $n_i< n_s$}\\
\nu_i  & \hbox{if $c-n_s+n_i-1-d\leqslant\nu_i< c$}\\
c  & \hbox{if $\nu_i<c-n_s+n_i-1-d$ and $n_i\geqslant n_s+d+1$}\\
c-n_s+n_i-1-d  & \hbox{if $\nu_i<c-n_s+n_i-1-d$ and $n_i< n_s+d+1$,}
\end{cases}
$$

We put $S_2(d)=\{i\mid 1\leqslant i\leqslant k, \mu_i<c-n_s+n_i-1-d, n_i< n_s+d+1\}$.

Then, we easily see  $S_2(d)= \{i\in\underline{\Psi}(\breve{\mu}^{(d)};c-n_s-d)\mid c>\breve{\mu}^{(d)}_i\}$ and
$$\breve{\mu}^{(d-1)}=\breve{\mu}^{(d)}+\sum_{i\in S_2(d)}f_i,$$
 for $1\leqslant d\leqslant c-n_s-b.$
From the definition of  $\breve{\mu}^{(d)}$, we see
$$\Phi(\breve{\mu}^{(d)};c-n_s-d)-S_2(d)=\{j\mid 1\leqslant j\leqslant k, c=\breve{\mu}^{(d)}_j, n_i\geqslant n_s+d+1\}.$$

For any $i\in S_2(d)$, we have
 $$c-n_s-d+n_i-1<c$$
So, we can apply Lemma 3.2.6 and we have
$$T_{\breve{\mu}^{(d)}}^{\breve{\mu}^{(d-1)}}\left(M_\T\left(\breve{\mu}^{(d)}\right)\right)\cong M_\T\left(\breve{\mu}^{(d-1)}\right).$$
Similarly, we have
$$T_{\breve{\nu}^{(d)}}^{\breve{\nu}^{(d-1)}}\left(M_\T\left(\breve{\nu}^{(d)}\right)\right)\cong M_\T\left(\breve{\nu}^{(d-1)}\right).$$
From Lemma 3.2.4, we have
$$T_{\breve{\mu}^{(d)}}^{\breve{\mu}^{(d-1)}}=T_{w\breve{\mu}^{(d)}}^{w\breve{\mu}^{(d-1)}}=T_{\breve{\nu}^{(d)}}^{\breve{\nu}^{(d-1)}}.$$
We have $\breve{\mu}^{(c-n_s-b)}=\widetilde{\mu}^{(0)}$ and  $\breve{\nu}^{(c-n_s-b)}=\widetilde{\nu}^{(0)}$.

Since translation functors are exact, we have
 $$M_\T\left(\breve{\nu}^{(0)}\right)\subseteq
 M_\T\left(\breve{\mu}^{(0)}\right)$$
 from $M_\T\left(\widetilde{\nu}^{(0)}\right)\subseteq
 M_\T\left(\widetilde{\mu}^{(0)}\right).$ 

\medskip

{\it Step 5}  \,\,\,\,

For $0\leqslant d\leqslant n_s-1$, we define $\bar{\mu}^{(d)}=[\bar{\mu}^{(d)}_1,...,\bar{\mu}^{(d)}_k]$ and  $\bar{\nu}^{(d)}=[\bar{\nu}^{(d)}_1,...,\bar{\nu}^{(d)}_k]$ as follows.

$$\bar{\mu}^{(d)}_i=
\begin{cases}
c  & \hbox{if  $n_i\geqslant n_s$}\\
c  & \hbox{if  $n_i= n_s$ and $\mu_i\geqslant c$}\\
c-1  & \hbox{if  $n_i= n_s$ and $\mu_i< c$}\\
c-1  & \hbox{if $d< n_i< n_s$}\\
c-1  & \hbox{if $n_i\leqslant d$ and $\mu_i\geqslant c$}\\
\mu_i  & \hbox{if   $n_i\leqslant d$ and  $c+n_i-2-d\leqslant\mu_i< c$}\\
c+n_i-2-d  & \hbox{if  $n_i\leqslant d$ and $c+n_i-2-d>\mu_i$,}
\end{cases}
$$
$$\bar{\nu}^{(d)}_i=
\begin{cases}
c  & \hbox{if  $n_i\geqslant n_s$}\\
c  & \hbox{if  $n_i= n_s$ and $\nu_i\geqslant c$}\\
c-1  & \hbox{if  $n_i= n_s$ and $\nu_i< c$}\\
c-1  & \hbox{if $d<n_i< n_s$}\\
c-1  & \hbox{if $n_i\leqslant d$ and $\nu_i\geqslant c$}\\
\nu_i  & \hbox{if   $n_i\leqslant d$ and  $c+n_i-2-d\leqslant\nu_i< c$}\\
c+n_i-2-d  & \hbox{if  $n_i\leqslant d$ and $c+n_i-2-d>\nu_i$.}
\end{cases}
$$

We put $S_3(d)=\{i\mid 1\leqslant i\leqslant k, \mu_i<c+n_i-2-d, n_i\leqslant d\}$.

Then, we easily see  $S_3(d)= \{i\in\underline{\Psi}(\bar{\mu}^{(d)};c-1-d)\mid n_i\leqslant d\}$ and
$$\bar{\mu}^{(d-1)}=\bar{\mu}^{(d)}+\sum_{i\in S_3(d)}f_i,$$
 for $1\leqslant d\leqslant n_s-1.$
From the definition of  $\bar{\mu}^{(d)}$, we see
\begin{multline*}
\Phi(\bar{\mu}^{(d)};c-1-d)-S_3(d)= \\ \{j\mid 1\leqslant j\leqslant k, c-1=
\bar{\mu}^{(d)}_j, n_i\geqslant d+1\}\cup\{j\mid 1\leqslant j\leqslant k, c=\bar{\mu}^{(d)}_j, n_i\geqslant n_s\}.
\end{multline*}
For any $i\in S_3(d)$, we have
 $$c+n_i-2-d<c-1.$$
So, we can apply Lemma 3.2.6 and we have
$$T_{\bar{\mu}^{(d)}}^{\bar{\mu}^{(d-1)}}\left(M_\T\left(\bar{\mu}^{(d)}\right)\right)\cong M_\T\left(\bar{\mu}^{(d-1)}\right).$$
Similarly, we have
$$T_{\bar{\nu}^{(d)}}^{\bar{\nu}^{(d-1)}}\left(M_\T\left(\bar{\nu}^{(d)}\right)\right)\cong M_\T\left(\bar{\nu}^{(d-1)}\right).$$
From Lemma 3.2.4, we have
$$T_{\bar{\mu}^{(d)}}^{\bar{\mu}^{(d-1)}}=T_{w\bar{\mu}^{(d)}}^{w\bar{\mu}^{(d-1)}}=T_{\bar{\nu}^{(d)}}^{\bar{\nu}^{(d-1)}}.$$

We assume that $n_s<\max\{n_i\mid 1\leqslant i\leqslant k\}$.
Then, we have $\bar{\mu}^{(n_s-1)}=\breve{\mu}^{(0)}$ and  $\bar{\nu}^{(n_s-1)}=\breve{\nu}^{(0)}$.
Since translation functors are exact, we have
 $$M_\T\left(\bar{\nu}^{(0)}\right)\subseteq
 M_\T\left(\bar{\mu}^{(0)}\right)$$
 from $M_\T\left(\breve{\nu}^{(0)}\right)\subseteq
 M_\T\left(\breve{\mu}^{(0)}\right).$ 

Next, we assume there is no  $1\leqslant i\leqslant k$ such that $n_s<n_i$.
In this case, we have  $\bar{\mu}^{(n_s-1)}=\widetilde{\mu}^{(0)}$ and  $\bar{\nu}^{(n_s-1)}=\widetilde{\nu}^{(0)}$.
Since translation functors are exact, we also have
 $$M_\T\left(\bar{\nu}^{(0)}\right)\subseteq
 M_\T\left(\bar{\mu}^{(0)}\right)$$
 from $M_\T\left(\widetilde{\nu}^{(0)}\right)\subseteq
 M_\T\left(\widetilde{\mu}^{(0)}\right).$ 

Let  $\bar{\mu}$ and  $\bar{\nu}$ be as in the statement of Lemma 3.1.3.
Since we immediately see $\bar{\mu}=\bar{\mu}^{(0)}$ and  $\bar{\nu}=\bar{\nu}^{(0)}$, anyway we have
 $$M_\T\left(\bar{\nu}\right)\subseteq
 M_\T\left(\bar{\mu}\right).$$
On the other hand, we have 
 $$M_\T\left(\bar{\nu}\right)\not\subseteq
 M_\T\left(\bar{\mu}\right)$$
from Lemma 3.1.3.
Therefore, we deduce a contradiction, as we desired.   
\,\,\,\, Q.E.D.

\subsection{An example}

In order to illustrate the above proof of Theorem 2.4.3, we consider the following example.
Let $\ggg=\gl(14,\cpx)$.
We put $k=6$ and $(n_1,n_2,n_3,n_4,n_5,n_6)=(4,1,2,1,2,4)$.
We consider the corresponding $\T=\Pi-\{\alpha_4, \alpha_5, \alpha_7,\alpha_8,\alpha_{10}\}$.
We also put 
\[ \mu=[4,3,-1,-4,2,-2],\,\,\,\,  \nu=[-2,-4,2,3,-1,4].\]
Applying the arguments in \S 3, we check that $M_\T(\nu)\not\subseteq M_\T(\mu)$.
We see that $c=2$ and $s=3$ (cf.\ Lemma 3.1.2).   So, $n_s=2$.
We have 
\[ \bar{\mu}=[2,1,1,1,2,2],\,\,\,\,  \bar{\nu}=[2,1,2,1,1,2], \]
and $M_\T(\bar{\nu})\not\subseteq M_\T(\bar{\mu})$  (cf.\ Lemma 3.1.3).

We assume $M_\T(\nu)\subseteq M_\T(\mu)$ and derive an impossible inclusion $M_\T(\bar{\nu})\subseteq M_\T(\bar{\mu})$.

(Step1) \,\, In this step, we make the entries greater than $c=2$ into $2$.
\[ \mu=\hat{\mu}^{(2)}=[\,\hbox{\fbox{$4$}},3,-1,-4,2,-2], \,\,  \nu=\hat{\nu}^{(2)}=[-2,-4,2,3,-1,\hbox{\fbox{$4$}}\, ],\]
\[ \hat{\mu}^{(1)}=[\,\hbox{\fbox{$3$}},\hbox{\fbox{$3$}},-1,-4,2,-2], \,\,  \hat{\nu}^{(1)}=[-2,-4,2,\hbox{\fbox{$3$}},-1,\hbox{\fbox{$3$}}\, ],\]
\[ \hat{\mu}^{(0)}=[2,2,-1,-4,2,-2], \,\,  \hat{\nu}^{(0)}=[-2,-4,2,2,-1,2],\]

(Step2) \,\,  For $1\leqslant i\leqslant 6$ such that $n_i<2$, if the $i$-th entry is $2$, we replace that by $1$.
\[ \hat{\mu}^{(0)}=[2,\hbox{\fbox{$2$}},-1,-4,2,-2], \,\,  \hat{\nu}^{(0)}=[-2,-4,2,\hbox{\fbox{$2$}},-1,2],\]
\[ \widehat{\mu}=[2,1,-1,-4,2,-2], \,\,  \widehat{\nu}=[-2,-4,2,1,-1,2].\]

In the remaining three steps, we take care of entries less than 2.

(Step3) \,\, We remark that $4=\max\{n_i\mid 1\leqslant i\leqslant 6\}$.  In this step,  if the $i$-th entry is smaller than $c-4+n_i$, we make that into $c-4+n_i$.
\[ \widehat{\mu}=\widetilde{\mu}^{(4)}=[2,1,-1,-4,2,\hbox{\fbox{$-2$}}], \,\,  \widehat{\nu}=\widetilde{\nu}^{(4)}=[\,\hbox{\fbox{$-2$}},-4,2,1,-1,2],\]
\[\widetilde{\mu}^{(3)}=[2,1,-1,\hbox{\fbox{$-4$}},2,\hbox{\fbox{$-1$}}], \,\,  \widetilde{\nu}^{(3)}=[\, \hbox{\fbox{$-1$}},\hbox{\fbox{$-4$}},2,1,-1,2],\]
\[\widetilde{\mu}^{(2)}=[2,1,-1,\hbox{\fbox{$-3$}},2,\hbox{\fbox{$0$}}], \,\,  \widetilde{\nu}^{(2)}=[\, \hbox{\fbox{$0$}},\hbox{\fbox{$-3$}},2,1,-1,2],\]
\[\widetilde{\mu}^{(1)}=[2,1,\hbox{\fbox{$-1$}},\hbox{\fbox{$-2$}},2,\hbox{\fbox{$1$}}], \,\,  \widetilde{\nu}^{(1)}=[\, \hbox{\fbox{$1$}},\hbox{\fbox{$-2$}},2,1,\hbox{\fbox{$-1$}},2],\]
\[\widetilde{\mu}^{(0)}=[2,1,0,-1,2,2], \,\,  \widetilde{\nu}^{(0)}=[2,-1,2,1,0,2].\]

(Step4) \,\,  In this step, we perform the following two procedures.
\begin{itemize}
\item  If the $i$-th entry is smaller than $c$ and $n_i>n_s$ , we make that into $c$.
\item If the $i$-th entry is smaller than $c-n_i+n_s$ and $n_i\leqslant n_s$ , we make that into  $c-n_i+n_s$.
\end{itemize}

$$\widetilde{\mu}^{(0)}=\breve{\mu}^{(1)}=[2,1,\hbox{\fbox{$0$}},\hbox{\fbox{$-1$}},2,2], \,\, \widetilde{\nu}^{(0)}=\breve{\nu}^{(1)}=[2,\hbox{\fbox{$-1$}},2,1,\hbox{\fbox{$0$}},2],$$
$$\breve{\mu}^{(0)}=[2,1,1,0,2,2], \,\, \breve{\nu}^{(0)}=[2,0,2,1,1,2].$$

(Step5) \,\,  In this step, if the $i$-th entry is smaller than $c-1$ and $n_i<n_s$, we make that into  $c-1$.
$$\breve{\mu}^{(0)}=\bar{\mu}^{(1)}=[2,1,1,\hbox{\fbox{$0$}},2,2], \,\, \breve{\nu}^{(0)}=\bar{\nu}^{(1)}=[2,\hbox{\fbox{$0$}},2,1,1,2],$$
$$\bar{\mu}=\bar{\mu}^{(0)}=[2,1,1,1,2,2], \,\, \bar{\nu}=\bar{\nu}^{(0)}=[2,1,2,1,1,2].$$
So, we have an impossible inclusion  $M_\T(\bar{\nu})\subseteq M_\T(\bar{\mu})$.

\end{document}